\documentclass[preprint,12pt]{article}

 \usepackage{graphics}
 \usepackage{graphicx}
 \usepackage{epsfig}
 \usepackage[all]{xy}

\usepackage{amssymb}
\usepackage{amsthm,amsmath}

\usepackage{latexsym, epsfig, psfrag,eepic,colordvi,bm}
\usepackage{graphics,color}
\usepackage{amsmath,amsfonts,amssymb,amscd}
\usepackage[all]{xy}
\usepackage{epstopdf}
\usepackage{cite}

\usepackage{amsthm,amsmath,amssymb}

\usepackage{graphicx}

\theoremstyle{plain}
\newtheorem{theorem}{Theorem}[section]
\newtheorem{lemma}[theorem]{Lemma}
\newtheorem{corollary}[theorem]{Corollary}
\newtheorem{proposition}[theorem]{Proposition}

\theoremstyle{definition}
\newtheorem{definition}[theorem]{Definition}

\theoremstyle{remark}

\newcommand{\bai}{\hspace{6pt}}

\begin{document}

		

\title{ Combinatorially refine a Zagier-Stanley result on products of permutations }

\author{Ricky X. F. Chen\\
	\small Biocomplexity Institute and Initiative, University of Virginia\\[-0.8ex]
	\small 995 Research Park Blvd, Charlottesville, VA 22911, USA\\[-0.8ex]
	\small\tt chenshu731@sina.com
}

\date{}
\maketitle


\begin{abstract}
In this paper, we enumerate the pairs of permutations that are long cycles and
whose product has a given cycle-type.
Our main result is a simple relation concerning the desired numbers for a 
few related cycle-types.
The relation refines a formula of the number of pairs of long cycles
whose product has $k$ cycles independently obtained by Zagier and Stanley 
relying on group characters, and was previously obtained by F\'{e}ray and Vassilieva by counting some colored permutations first and then relying on some algebraic computations in the ring of symmetric functions.
Our approach here is simpler and combinatorial.

  

\end{abstract}



\section{Introduction}

Let $\mathfrak{S}_n$ denote the symmetric group on $[n]=\{1,2,\ldots, n\}$.
We shall use the following two representations of a permutation $\pi\in \mathfrak{S}_n$:\\
\emph{two-line form:} the top line lists all elements in $[n]$, following the natural order.
The bottom line lists the corresponding images of elements on the top line, i.e.,
\begin{eqnarray*}
	\pi=\left(\begin{array}{ccccccc}
		1&2& 3&\cdots &n-2&{n-1}&n\\
		\pi(1)&\pi(2)&\pi(3)&\cdots &\pi({n-2}) &\pi({n-1})&\pi(n)
	\end{array}\right).
\end{eqnarray*}
\emph{cycle form:} a permutation $\pi$ is decomposed into disjoint cycles.
The set consisting of the lengths of these disjoint cycles is called the cycle-type of $\pi$. 
We can encode this set as an integer partition of $n$.
An integer partition $\lambda$ of $n$, denoted by $\lambda \vdash n$,
can be represented by a non-increasing integer sequence $\lambda=\lambda_1
\lambda_2\cdots$, where $\sum_i \lambda_i=n$, or as $1^{m_1}2^{m_2}\cdots n^{m_n}$, where
we have $m_i(\lambda)$ of part $i$ and $\sum_i i m_i =n$.
 A cycle of length $k$ is called a $k$-cycle.
In addition, we denote the number of permutations of cycle-type $\lambda$ by $z_{\lambda}$.
It is well known that $z_{\lambda}=\frac{n!}{\prod_i i^{m_i}  m_i !}$ if $\lambda=1^{m_1} 2^{m_2} \cdots n^{m_n}$.	
We also denote the length of $\lambda$, i.e.,~the number of positive parts in $\lambda$, by $\ell(\lambda)$.

Factorizations of permutations or products of permutations
have been extensively studied in different contexts.
Most of related results in the field
 rely either partially or totally on a character theoretic approach (e.g.,~\cite{jack, stan, zag}). 
It is generally hard to obtain explicit and simple counting formulas.
However, when one of the involved permutations is a long cycle,
we may obtain some explicit formulas, see~\cite{chr-1,chapuy,G-S,walsh-lehman} and references therein.
In particular,
Zagier~\cite{zag} and Stanley~\cite{stan3}
have independently obtained the following result: the number of $n$-cycles $s$ such that the product $(1\bai 2\bai \cdots \bai n)\, s$ has $k$ cycles is given by the surprisingly simple formula $\frac{2}{n(n+1)} C(n+1,k)$,
where $C(n,k)$ stands for the signless Stirling number of the first kind, i.e., the number of permutations on $[n]$ with $k$ cycles.
Combinatorial proofs of this result can be found in~\cite{cms,fv,chr-1}.
Recently, the author~\cite{chen3} also obtained some analogues of the above Zagier-Stanley formula in
the context of studying separation probabilities~\cite{bdms, stan2} of products of permutations.
For instance, it was proved that the number of pairs of $n$-cycles whose product has $k$ cycles and separates the elements in $[m]$ is given by 
$$
\frac{2 (n-1)!}{(n+m)(n+1-m)}  C_m(n+1,k),
$$
where $C_m(n,k)$ is the number of permutations on $[n]$
with $k$ cycles and the elements in $[m]$ separated, i.e., an analogue of $C(n,k)$.
This analogue was particularly used to answer a call of Stanley~\cite{stan2} for simple combinatorial proofs for the probability of separating $m$ elements due to Du and Stanley~\cite{stan2}.

In this paper, we enumerate the pairs of long cycles whose product has a given cycle-type,
refining the Zagier-Stanley result.
Specifically, we obtain the theorem described below.

Let $\lambda=1^{m_1} 2^{m_2} \cdots n^{m_n}\vdash n+1$.
For $i>0$ (and $m_{i+1}\neq 0$), denote by $\lambda^{\downarrow (i+1)}$
the partition $\mu=1^{m_1} \cdots i^{m_i+1} (i+1)^{m_{i+1}-1} \cdots n^{m_n} \vdash n$, i.e.,~changing
an $i+1$ part to an $i$ part. 
Let $p^{(n)}_{\mu}$ denote the number of pairs of $n$-cycles whose product
has a cycle-type $\mu$.

\begin{theorem}\label{thm:main}
Suppose $m$ and $n$ have the same parity. Then for any partition $\lambda \vdash n+1$ of length $m$, we have
\begin{align}\label{eq:thm-main}
\frac{n+1}{2} \sum_{\mu=\lambda^{\downarrow (i+1)}, \; i>0} i m_i(\mu) p^{(n)}_{\mu} = (n-1)! z_{\lambda} \; .
\end{align}
\end{theorem}

Summing over all possible partitions of length $m$ in eq.~\eqref{eq:thm-main}
will give us the Zagier-Stanley result.
An equivalent statement of Theorem~\ref{thm:main}
was previously obtained in F\'{e}ray and Vassilieva~\cite{fv}, 
by 
counting some colored permutations first and then by some algebraic computations in the ring of symmetric functions.
Our approach here is simpler and totally combinatorial,
and is based on extending
the plane permutation framework which was first introduced by the author and Reidys~\cite{chr-1} and has proven to be effective in dealing with hypermaps, graph embeddings and even the genome rearrangement problems
involving transpositions, block-interchanges and reversals~\cite{chr-1, chr-2, chen3}.

\section{Refining the Zagier-Stanley result}

We begin with a review of some notation and results on plane permutations in~\cite{chr-1}.
%%%
%%%%%%%%%%%%%%%%%%%%%%%%%%%%%%%%%%%%%%%%%%%%%%%%%%%%%%%%%%%%%%%%%%%%%%%%%%%%%%%%%%%%%
%%%
%%%
%%%%%%%%%%%%%%%%%%%%%%%%%%%%%%%%%%%%%%%%%%%%%%%%%%%%%%%%%%%%%%%%%%%%%%%%%%%%%%%%%%%%%
%%%
\begin{definition}
	A \emph{plane permutation} on $[n]$ is a pair $\mathfrak{p}=(s,\pi)$ where $s=(s_i)_{i=0}^{n-1}$
	is an $n$-cycle and $\pi$ is an arbitrary permutation on $[n]$.
	Given $s=(s_0~s_1~\cdots ~s_{n-1})$,
a plane permutation $\mathfrak{p}=(s,\pi)$ is represented by a two-row array:
\begin{equation}
\mathfrak{p}=\left(\begin{array}{ccccc}
s_0&s_1&\cdots &s_{n-2}&s_{n-1}\\
\pi(s_0)&\pi(s_1)&\cdots &\pi(s_{n-2}) &\pi(s_{n-1})
\end{array}\right).
\end{equation}
The permutation $D_{\mathfrak{p}}$ induced by the diagonal-pairs (cyclically), i.e.,~$D_{\mathfrak{p}}(\pi(s_{i-1}))=s_i$ for $0<i< n$, and
$D_{\mathfrak{p}}(\pi(s_{n-1}))=s_0$, is called the \emph{diagonal} of $\mathfrak{p}$.

\end{definition}\label{2def1}
%%%
%%%%%%%%%%%%%%%%%%%%%%%%%%%%%%%%%%%%%%%%%%%%%%%%%%%%%%%%%%%%%%%%%%%%%%%%%%%%%%%%%%%%%
%%%

We sometimes refer to $s,\, \pi, \, D_{\mathfrak{p}}$ respectively as the upper horizontal, the vertical and the diagonal.
Obviously, we have
$D_{\mathfrak{p}}=
s \pi^{-1}$.
It should be pointed out that, although as a cyclic permutation, there is no absolute left-right order for the elements in $s$,
in this paper, we generally assume there is a left-right order, with the leftmost element being $s_0$.

In a permutation $\pi$ on $[n]$, $i$ is called an \emph{exceedance} if $i<\pi(i)$ following the natural order and an \emph{anti-exceedance} otherwise.
Note that $s$ induces a linear order $<_s$,
where $a<_s b$ if $a$ appears before $b$ in $s$ from left to right (with the leftmost element $s_0$).
Without loss of generality, we always assume $s_0=1$ unless explicitly stated otherwise.
These concepts then can be generalized for plane permutations as follows:
%%%
%%%%%%%%%%%%%%%%%%%%%%%%%%%%%%%%%%%%%%%%%%%%%%%%%%%%%%%%%%%%%%%%%%%%%%%%%%%%%%%%%%%%%
%%%
\begin{definition}\label{2def2}
	For a plane permutation $\mathfrak{p}=(s,\pi)$, an element $s_i$ is called an
	\emph{exceedance} of $\mathfrak{p}$ if $s_i<_s \pi(s_i)$, and an \emph{anti-exceedance} if $s_i\ge_s \pi(s_i)$.
\end{definition}
%%%
%%%%%%%%%%%%%%%%%%%%%%%%%%%%%%%%%%%%%%%%%%%%%%%%%%%%%%%%%%%%%%%%%%%%%%%%%%%%%%%%%%%%%
%%%

In the following, we mean by ``the cycles of $\mathfrak{p}=(s,\pi)$'' the cycles of $\pi$ and
any comparison of elements in $s,~\pi$ and $D_{\mathfrak{p}}$ references the linear order $<_s$.

Obviously, each $\mathfrak{p}$-cycle contains at least one anti-exceedance as it contains
a minimum, $s_i$, for which $\pi^{-1}(s_i)$ is an anti-exceedance. We call these trivial anti-exceedances
and refer to a \emph{non-trivial anti-exceedance} as an NTAE. Furthermore, in any cycle of length
greater than one, its minimum is always an exceedance.

Let $\mu,~\lambda$ be two integer partitions of $n$.
We denote $\mu\rhd_{k} \lambda$ if $\mu$ can be obtained from $\lambda$ by splitting one part into $k$ parts,
or equivalently, $\lambda$ from $\mu$ by merging $k$ parts into one part.
Let $\kappa_{\mu,\lambda}$ be the number of different ways of merging $k$ parts of $\mu$
in order to obtain $\lambda$ provided that $\mu \rhd_k \lambda$.
Note that we differentiate two parts of $\mu$ even if the two parts are of the same value.
For example, for $\mu=1^2 2^2$ and $\lambda=1^1 2^1 3^1$, we have
$\kappa_{\mu, \lambda}=4$.

Let $U_{\lambda}^{\eta}$ denote the set of plane permutations on $[n]$
where the diagonal is of cycle-type $\eta$ and the vertical is of cycle-type $\lambda$. 
We always assume $\ell(\lambda)+\ell(\eta)$ has the same parity as $n+1$.
Otherwise we know $U_{\lambda}^{\eta}=\varnothing$.
We denote $p^{\eta}_{\lambda}= |U_{\lambda}^{\eta}|$.	
In Chen and Reidys~\cite{chr-1}, 
 while studying the transposition action on the diagonal of plane permutations (i.e.,~transposing
two adjacent diagonal blocks where a diagonal block is a set of consecutive diagonal-pairs)
and exceedances, motivated by the work~\cite{chapuy},
 the following proposition has been proved.

\begin{proposition}[\cite{chr-1}]
Let $\lambda, \eta \vdash n$ and $p^{\eta}_{\lambda,a}$ be the number of $\mathfrak{p} \in U^{\eta}_{\lambda}$ such that $\mathfrak{p}$
has $a$ exceedances. Then
we have
\begin{align}
\sum_{a\geq 0} \bigl(n-\ell(\lambda)-a \bigr) p^{\eta}_{\lambda,a} &=\sum_{ \mu \rhd_{2i+1} \lambda, \; i>0} \kappa_{\mu,\lambda} p^{\eta}_{\mu} \label{eq:gen}\; ,\\
\bigl( n+1-\ell(\lambda) \bigr) p^{(n)}_{\lambda} &=\sum_{\mu \rhd_{2i+1} \lambda, \; i>0} \kappa_{\mu,\lambda} p^{(n)}_{\mu} + (n-1)! z_{\lambda} \label{eq:long}\; .
\end{align}
\end{proposition}

Our new contribution starts from here.
In eq.~\eqref{eq:gen}, if we sum over all $\eta \vdash n$, we realize that
\begin{align}\label{cor:exc}
\sum_{\eta\vdash n}  \sum_{a\geq 0} a p^{\eta}_{\lambda,a} =\bigl( n-\ell(\lambda)\bigr) (n-1)! z_{\lambda} - \sum_{ \mu \rhd_{2i+1} \lambda, \; i>0} \kappa_{\mu,\lambda} (n-1)! z_{\mu} \; .
\end{align}

Note that the left-hand side of the above equation can be interpreted as the total number of exceedances of
$\mathfrak{p} \in \bigcup_{\eta \vdash n} U^{\eta}_{\lambda}$.
However, by directly counting the exceedances, we have

\begin{lemma}\label{lem:exc}
The total number of exceedances
\begin{align}\label{eq:exc}
\sum_{\eta \vdash n}  \sum_{a\geq 0} a p^{\eta}_{\lambda, a} =\frac{n-m_1(\lambda)}{2} (n-1)! z_{\lambda} \; .
\end{align}
\end{lemma}
\proof 
It should not be hard to observe that in $\bigcup_{\eta \vdash n} U^{\eta}_{\lambda}$, for any upper horizontal $s$, the total number of exceedances
of the plane permutations with the upper horizontal $s$ is the same as the total number of 
exceedances of the plane permutations with the upper horizontal $(1\bai 2\bai \cdots \bai n)$.
The latter is really just counting the total number of exceedances of the (conventional) permutations
of cycle-type $\lambda$.
Note that if a permutation $\pi$ of cycle-type $\lambda$ has $a$ exceedances, then its inverse $\pi^{-1}$
is of cycle-type $\lambda$ with $n-m_1(\lambda)-a$ exceedances. Because if an element $x$ that is not a fixed point
is an exceedance of $\pi$, $\pi(x)$ is an exceedance of $\pi^{-1}$ and not a fixed point.
Thus, for each such pair, they have on average $\frac{n-m_1(\lambda)}{2} $ exceedances,
completing the proof. \qed

Combining eq.~\eqref{cor:exc} and eq.~\eqref{eq:exc}, we have

\begin{proposition} For any $\lambda \vdash n+1$, the following is true
\begin{align}\label{eq:baserecur}
\bigl( n+1-\ell(\lambda) \bigr) z_{\lambda} 
=\sum_{\mu \rhd_{2i+1} \lambda, \; i>0} \kappa_{\mu,\lambda} z_{\mu} + \frac{ z_{\lambda} }{2} \sum_{i>0} (i+1) m_{i+1} (\lambda) \; .
\end{align}
\end{proposition}
\proof Just note that $n+1-m_1(\lambda)=\sum_{i>0} (i+1) m_{i+1}(\lambda)$. \qed

As a consequence of eq.~\eqref{eq:long}, we have the following corollary.

\begin{corollary} For any $\lambda \vdash n+1$ and $i>0$, we obtain
\begin{multline}\label{eq:downarrow}
\bigl( n+1-\ell(\lambda) \bigr) \frac{n+1}{2} i m_i (\lambda^{\downarrow (i+1)}) p^{(n)}_{\lambda^{\downarrow (i+1)}} \\
=\sum_{j>0, \atop \mu \rhd_{2j+1} \lambda^{\downarrow (i+1)}} \kappa_{\mu,\lambda^{\downarrow (i+1)}}  \frac{n+1}{2} i m_i (\lambda^{\downarrow (i+1)}) p^{(n)}_{\mu} 
 +\frac{(i+1) m_{i+1} (\lambda )}{2} (n-1)! z_{\lambda} \; .
\end{multline}
\end{corollary}
\proof Based on eq.~\eqref{eq:long}, we first have
$$
\bigl( n+1-\ell(\lambda) \bigr) p^{(n)}_{\lambda^{\downarrow (i+1)}} =\sum_{\mu \rhd_{2j+1} \lambda^{\downarrow (i+1)}, \; j>0} \kappa_{\mu,\lambda^{\downarrow (i+1)}} p^{(n)}_{\mu} + (n-1)! z_{\lambda^{\downarrow (i+1)}}.
$$
Next, we observe 
$
z_{\lambda^{\downarrow (i+1)}} = \frac{(i+1) m_{i+1} (\lambda)}{(n+1) i m_i (\lambda^{\downarrow (i+1)})} z_{\lambda} ,
$
and the proof follows. \qed

In order to proceed, we need the following key lemma.

\begin{lemma} For any $\lambda \vdash n+1$, it holds that
\begin{align}\label{eq:keylem}
\sum_{i>0} i m_i (\lambda^{\downarrow (i+1)}) \sum_{j>0, \atop \mu \rhd_{2j+1} \lambda^{\downarrow (i+1)}} \kappa_{\mu,\lambda^{\downarrow (i+1)}}    p^{(n)}_{\mu} 
= \sum_{j>0, \atop \mu \rhd_{2j+1} \lambda} \kappa_{\mu,\lambda}  \sum_{i>0} i m_i (\mu^{\downarrow (i+1)}) p^{(n)}_{\mu^{\downarrow (i+1)}}.
\end{align}
\end{lemma}
\proof Note that both sides are eventually sums where each summand is indexed by a partition $\mu$ of $n$.
So it suffices to show the coefficients of the $\mu$-summands on both sides agree.
Suppose $\mu=\mu_1 \mu_2 \cdots$. We identify $\mu$ as the permutation
$\mu=(1\bai 2\bai \cdots \bai \mu_1)(\mu_1+1 \bai \mu_1+2 \bai \cdots \bai \mu_1+\mu_2)\cdots$.
Note that the number $\kappa_{\mu, \beta}$ for $\mu \rhd_{2k+1} \beta$ can be interpreted as 
the number of permutations of cycle-type $\beta$ that can be obtained by concatenating $2k+1$ cycles
of the permutation $\mu$ into one single cycle according to the cycle lengths (decreasingly) and the minimum
elements of the cycles (increasingly).
Next, we shall show that the considered coefficients counting the same subset of permutations of cycle-type $\lambda$ on $[n+1]$
associated with the permutation $\mu$.

Suppose $\ell(\lambda)= \ell(\mu)-2j$. Let us start with the left-hand side and describe the subset $A$ of permutations associated with $\mu$.
A permutation $\gamma \in A$ if $\gamma$ gives us the permutation $\mu$ following the procedure:
(i) erasing $n+1$ from the cycle of $\gamma$ containing it to obtain $\gamma'$; (ii) splitting one cycle of $\gamma'$
into $2j+1$ cycles.
Now suppose $\mu \rhd_{2j+1} \lambda^{\downarrow (i+1)}$.
We claim that there are $i m_i (\lambda^{\downarrow (i+1)}) \kappa_{\mu,\lambda^{\downarrow (i+1)}} $
associated permutations contained in $A$.
This can be seen from the other way around:
merging $2j+1$ cycles of $\mu$ to obtain a permutation $\gamma'$ of cycle-type
$ \lambda^{\downarrow (i+1)}$ in $\kappa_{\mu,\lambda^{\downarrow (i+1)}}$
different ways; and next inserting $n+1$ into a cycle of length $i$ in $\gamma'$
to obtain a permutation $\gamma$ of cycle-type $\lambda$ in $i m_i (\lambda^{\downarrow (i+1)})$
different ways. Denote the subset of permutations by $A_i$. 
Next it is obvious that for $i\neq j$, $A_i \bigcap A_j =\varnothing$, because
$n+1$ is in cycles of different lengths. Therefore, $A=\bigcup_{i>0} A_i$ and
the coefficient of the $\mu$-term on the left-hand side is $|A|$.

Denote $B$ the subset of associated permutations on the right-hand side,
where a permutation $\gamma \in B$ if $\gamma$ gives us the permutation $\mu$ following the procedure:
(1) splitting one cycle of $\gamma$
into $2j+1$ cycles to obtain $\gamma'$; and
(2) erasing $n+1$ from the cycle of $\gamma'$ containing it.
Now suppose $\mu=\mu'^{\downarrow (i+1)}$ and $ \mu' \rhd_{2j+1} \lambda$.
Analogously, we can conclude there are $i m_i (\mu'^{\downarrow (i+1)}) \kappa_{\mu', \lambda}$
different associated permutations. 
Denote the subset of permutations by $B_i$. 
We also have $B_i \bigcap B_j =\varnothing$. Hence,
the coefficient of the $\mu$-term on the right-hand side is $|B|$.
Finally, we can easily check that actually $A=B$. Therefore,
the lemma follows. \qed

For any $\lambda \vdash n+1$, we denote $\mathcal{T}_{\lambda}= \sum_{i>0} \frac{n+1}{2}i m_i(\lambda^{\downarrow (i+1)}) p^{(n)}_{\lambda^{\downarrow (i+1)}} $. Based on eq.~\eqref{eq:downarrow} and eq.~\eqref{eq:keylem},
we obtain:
\begin{align}\label{eq:recur}
\bigl(n+1-\ell(\lambda) \bigr) \mathcal{T}_{\lambda}=\sum_{\mu \rhd_{2j+1} \lambda} \kappa_{\mu, \lambda} \mathcal{T}_{\mu}
+\frac{(n-1)! z_{\lambda}}{2} \sum_{i>0} (i+1) m_{i+1} (\lambda) \; .
\end{align}
Now we are ready to prove our main theorem.

{\bf Proof of Theorem~\ref{thm:main}}. Based on eq.~\eqref{eq:baserecur} and eq.~\eqref{eq:recur}, we observe that both sides of the equality
in the theorem satisfy the same recurrence.
Then, it suffices to compare the respective initial conditions.
Note that the initial cases correspond to the cases that $\lambda=1^a 2^b$ and $a+2b=n+1$
for some $a\geq 0$ and $b\geq 0$.
In~\cite{chr-1}, it is proved that the number of factorizations of 
a fixed permutation of cycle-type $\lambda$ into two long cycles is $\frac{n!}{n+2-a-b}$.
Then we can compute
\begin{align*}
\frac{n+1}{2} \sum_{\mu=\lambda^{\downarrow (i+1)}, \; i>0} i m_i(\mu)  p^{(n)}_{\mu}
&=\frac{n+1}{2}  (a+1) p^{(n)}_{1^{a+1} 2^{b-1}}\\
&=\frac{n+1}{2}  (a+1) \frac{z_{1^{a+1} 2^{b-1}}}{(n-1)!}\frac{(n-1)! (n-1)!}{n+1-(a+1)-(b-1)}\\
&=\frac{(n-1)! (n+1)!}{1^a2^b a! b!}= (n-1)! z_{\lambda}.
\end{align*}
This completes the proof. \qed

\section*{Acknowledgments}

I would like to thank Christian Reidys for encouragements and support.

				%%%%%%%%%%%%%%%%%%%%%%%%%%%%%%%%%%%%%%%%%%%%%%%%%%%%%%%%%%%%%%%%%%%%%%%%%%%%%%%%%%%%%%%%%%%%%%%%%%%%%%%%%%%%%%%%%

\end{document}